\documentclass[11pt]{article}

\usepackage[utf8]{inputenc}
\usepackage{amsmath}
\usepackage{amsfonts,amssymb,latexsym,multicol}
\usepackage[francais]{babel}
\usepackage[T1]{fontenc}

\newcounter{fig}
\newtheorem{prop}{Proposition}

\newcommand{\expli}[1]{\quad\text{\footnotesize (#1)}}

\newcommand{\eps}{\varepsilon}
\newcommand{\fhi}{\varphi}
\newcommand{\ioe}{\leqslant}
\newcommand{\soe}{\geqslant}
\renewcommand{\le}{\leqslant}
\renewcommand{\ge}{\geqslant}
\newcommand{\vers}{\rightarrow}

\newcommand{\demi}{{\frac{1}{2}}}

\newcommand{\ci}{{\rm ci}}
\newcommand{\Ci}{{\rm Ci}}
\newcommand{\si}{{\rm si}}
\newcommand{\Si}{{\rm Si}}
\newcommand{\vp}{{\rm v.p.}}

\newcommand{\Fcal}{{\mathcal F}}

\newcommand{\Lcal}{{\mathcal L}}

\newcommand{\Wcal}{{\mathcal W}}

\newcommand{\Fgot}{{\mathfrak F}}

\newcommand{\Vgot}{{\mathfrak V}}

\newcommand{\Int}{{\mathbb Z}}

\newcommand{\Real}{{\mathbb R}}
\newcommand{\Com}{{\mathbb C}}

\newcommand{\fin}{\hfill$\Box$}
\newcommand{\dem}{\noindent {\bf D\'emonstration\ }}
\newcommand{\fine}{\tag*{\mbox{$\Box$}}}

\providecommand{\bysame}{\leavevmode ---\ }
\providecommand{\og}{``}
\providecommand{\fg}{''}
\providecommand{\smfandname}{et}
\providecommand{\smfedsname}{\'eds.}
\providecommand{\smfedname}{\'ed.}
\providecommand{\smfmastersthesisname}{M\'emoire}
\providecommand{\smfphdthesisname}{Th\`ese}

\title{Sur une équation fonctionnelle approchée due à J. R. Wilton}

\author{Michel Balazard et Bruno Martin}

\begin{document}
\maketitle

\begin{flushright}
\textit{Pour Micha Tsfasman et Serge Vladut, très amicalement.}
\end{flushright}

\begin{center}
  {\sc Abstract}
\end{center}
\begin{quote}
{\footnotesize We give a new proof of an approximate functional equation, due to J. R. Wilton, for a trigonometric sum involving the divisor function. This allows us to improve on Wilton's error term and to give an explicit formula for an unspecified function involved in the functional equation.
}
\end{quote}

\begin{center}
  {\sc Keywords}
\end{center}
\begin{quote}
{\footnotesize Approximate functional equations, divisor function, Mellin transform, Hilbert transform\\MSC classification : 11N37 (11L07)}
\end{quote}



\section{Introduction}

 Dans l'article \cite{Wilton}, publié en 1933 et dévolu entre autres à l'étude des séries trigonométriques
\begin{align}
\psi_1(x)&=\sum_{n\soe 1}\frac{\tau(n)}{n}\cos 2\pi nx\label{t105}\\
\psi_2(x)&=\sum_{n\soe 1}\frac{\tau(n)}{n}\sin 2\pi nx\label{t106}
\end{align}
o\`u $\tau(n)$ d\'esigne le nombre de diviseurs du nombre entier naturel $n$, Wilton démontre  pour la somme partielle 
\[
\psi(x,v)=\sum_{n\ioe v}\frac{\tau(n)}{n}e^{2\pi i nx}
\]
l'équation fonctionnelle approchée suivante. 
\begin{prop}[\cite{Wilton}, (2.2), Theorem 2, p. 223]\label{t103}
Soit $K>0$. Pour $0<x\ioe 1$, $v\soe 1$ et $x^2v > K$, on a
\begin{equation}\label{t0}
\psi(x,v) - x\overline{\psi(1/x,x^2v)}=\Fgot(x)+\demi\log^2\frac 1x+\Big(\gamma-\log 2\pi +\demi \pi i\Big)\log \frac 1x+O\big((x^2v)^{-1/5}\big),
\end{equation}
où $\gamma$ désigne la constante d'Euler, et $\Fgot$ une fonction continue sur le segment $[0,1]$. La constante implicite dans le $O$ ne dépend que de $K$.
\end{prop}

\smallskip

Ainsi la somme partielle $\psi(x,v)$, qui est de période $1$ en la variable $x$, a également un comportement simple sous l'effet du changement de variables $(x,v)\mapsto (1/x,x^2v)$. Une telle situation est propice à la description du comportement de cette somme partielle en termes du développement du nombre réel $x$ en fraction continue. Cette voie, ouverte voici un siècle par Hardy et Littlewood pour l'étude des séries trigonométriques associées aux fonctions thêta elliptiques (cf. \cite{HL}), est suivie par Wilton dans \cite{Wilton} pour aboutir à des critères nécessaires et suffisants de convergence pour $\psi_1(x)$ et $\psi_2(x)$.

Wilton n'explicite pas la fonction $\Fgot(x)$ (ce n'est pas nécessaire pour l'étude de la convergence des séries \eqref{t105} et \eqref{t106}) et indique que le terme d'erreur $O\big((x^2v)^{-1/5}\big)$ {\og can easily be sharpened\fg}.  Le but du présent travail est triple : proposer une démonstration de \eqref{t0} différente de celle de Wilton, expliciter la fonction $\Fgot$, et améliorer le terme d'erreur. 

La méthode de Hardy et Littlewood d'étude des sommes partielles d'une série oscillante \textit{via} une équation fonctionnelle approchée connaît depuis quelques années un regain d'intérêt (voir par exemple les articles de Rivoal \cite{rivoal-2012}, Rivoal et Roques \cite{rivoal-roques-2013}, Rivoal et Seuret \cite{rivoal-seuret-2014}). \'Etablir de telles équations fonctionnelles peut s'avérer particulièrement délicat, et la méthode que nous développons ici vient compléter la palette de techniques disponibles.


\medskip 
  
Afin d'énoncer notre résultat, nous rappelons la définition de la \emph{fonction d'autocorr\'elation multiplicative de la
fonction {\og partie fractionnaire\fg}}, introduite  par B\'aez-Duarte {\it et al.}
(cf. \cite{baez-duarte-all}) dans le contexte de l'\'etude du crit\`ere de Nyman pour
l'hypoth\`ese de Riemann. En désignant par $\lfloor t\rfloor$ la partie entière du nombre réel $t$ et par $\{t\}$ sa partie fractionnaire, égale à $t-\lfloor t\rfloor$, nous posons pour $x \soe 0$ :
\begin{equation*}
 A(x)= \int_0^\infty \{t \} \{x t \}\frac{dt}{t^2}. 
\end{equation*}
Rappelons que $A$ est continue sur $[0,\infty[$ et vérifie les relations asymptotiques
\begin{align}
A(x) &\sim \demi \log x \quad (x \vers \infty),\label{t80}\\
A(x) &\sim \demi x \log (1/x) \quad (x \vers 0),\label{t81}
\end{align}
où \eqref{t81} résulte de \eqref{t80} et de l'identité $A(x) =x A(1/x)$.

Nous définissons maintenant une fonction $\Fgot \, :]0,\infty[\vers \Com$ par les formules
\begin{align}
\Re\Fgot(x) &=2x\cdot\vp\int_0^{\infty}A(t)\frac{x^2+t^2}{(t-x)(t+x)^2}\, \frac{dt}{t}-x\left(\demi \log^2 x+(\log 2\pi -\gamma)\log \frac 1x -c_0\right)-c_0\label{t111}\\
\Im\Fgot(x)&= \pi \Big(A(x)+\frac{x}{2}\log x-\frac{x +1}{2}(\log 2\pi -\gamma) \Big)\label{t107}
\end{align}
avec
$$
c_0=\frac{\pi^2}{24}-\demi\log^2 2\pi+\demi\gamma^2+\gamma\log 2\pi+2\gamma_1,
$$
$\gamma_1$ désignant la constante de Stieltjes d'indice $1$, soit 
$\displaystyle{ \gamma_1= \lim_{n\to \infty} \Big(\sum_{k=1}^n \frac{\log k}{k}-\demi \log^2 n }\Big)$.

Rappelons que $\vp$ signifie {\og valeur principale\fg} (voir \S\ref{t115}).

\smallskip

Nous obtenons le résultat suivant. 
\begin{prop}\label{prop:eq-approchee-psi1}
La fonction $\Fgot$ est continue sur $]0,\infty[$ et se prolonge par continuité en $0$. De plus, si $K_1$ et $K_2$ sont des nombres positifs arbitraires, on a pour $0<x\ioe K_1$, $v>0$ et $x^2v \soe K_2$, 
\begin{multline}\label{t74}
\psi(x,v)-x\overline{\psi(1/x,x^2v)}=\Fgot(x)+\demi\log^2\frac 1x+(\gamma-\log 2\pi +\demi \pi i)\log \frac 1x+O\left(\frac{\log^2(2+x^2v)}{(x^2v)^{1/2}} \right).
\end{multline}
 La constante implicite dans le $O$ ne dépend que de $K_1$ et $K_2$.
\end{prop}

\medskip

La démonstration de Wilton de \eqref{t0} repose sur la formule sommatoire obtenue en 1904 par Voronoï (cf. \cite{voronoi1}, p. 209-210) :
\begin{equation}\label{t71}
\sum_n\tau(n)f(n)=\int_0^{\infty}(\log t+2\gamma)f(t)dt+\sum_n\tau(n)\Vgot f(n)
\end{equation}
où nous avons pris, pour simplifier, la fonction $f$ indéfiniment dérivable et à support compact dans $]0,\infty[$, et où
$$
\Vgot f(y)= 2\pi\int_0^{\infty}M_0(4\pi\sqrt{xy})f(x)dx,
$$
la fonction $M_0$ s'exprimant en termes de fonctions de Bessel (c'est la notation de Wilton ; Voronoï utilise la fonction $\alpha(x)=2\pi M_0(4\pi\sqrt{x})$).
L'année précédente, en 1903, Voronoï avait démontré la majoration  
\begin{equation}\label{t127}
\Delta(x)=\sum_{n\le x} \tau(n)-x(\log x+2\gamma-1)  \ll x^{1/3} \log x \quad (x\soe 2),
\end{equation}
grâce à une généralisation de la méthode de l'hyperbole de Dirichlet (cf. \cite{34.0231.03}). Une variante convenable de la relation \eqref{t71} permet de retrouver cette estimation (cf. \cite{chandrasekharan-70}, chapter VIII, \S5), qui peut également être démontrée grâce à la méthode de van der Corput (cf \cite{MR1366197}, théorème 6.11).

La déduction par Wilton de la  proposition \ref{t103} à partir de la formule sommatoire de Voronoï n'est pas immédiate ; elle comporte des calculs qui ne sont guère plus rapides que ceux de la démonstration que nous proposons. Celle-ci s'appuie néanmoins, comme la formule de Voronoï, sur l'équation fonctionnelle de la fonction $\zeta$ de Riemann, qui est le fait mathématique essentiel d'où découle, \emph{in fine}, l'équation fonctionnelle approchée de Wilton.

\smallskip

Notre démarche peut être résumée de la fa\c{c}on suivante. Dans un premier temps, nous démontrons au \S\ref{t126} l'existence d'une fonction $\Fgot$, continue sur $]0,\infty[$ et se prolongeant par continuité en $0$, telle que \eqref{t74} soit vérifiée. Pour cela, par une intégration par parties, nous ramenons l'étude de la somme partielle $\psi(x,v)$ à celle de l'intégrale
\begin{equation}\label{u136}
I(x,v)=2\int_0^v \frac{\Delta(t)}{t}e^{2\pi itx}dt.
\end{equation}
Nous transformons ensuite $I(x,v)$ en intégrale sur la droite $\Re s =\demi$ du plan complexe (dite droite \emph{critique}) grâce au théorème de Plancherel pour la transformation de Mellin, donnée par
\begin{equation}\label{t123}
Mf(s)=\int_0^{\infty}f(t)t^{s-1}dt.
\end{equation}
Le résultat est le suivant :
\begin{equation}\label{u5}
I(x,v)=\int_{\sigma=1/2}\frac{\zeta^2(s)}{s}\lambda(s,2\pi xv)(2\pi x)^{-s}\,\frac{d\tau}{\pi}.
\end{equation}
Dans cette relation interviennent d'une part la fonction $\zeta$ de Riemann, dont nous utiliserons des propriétés classiques rappelées au \S\ref{t75}, et d'autre part la fonction 
\begin{equation}
\label{t99}
\lambda(s,v) =\int_0^ve^{it}t^{s-1}dt.
\end{equation}
Cette dernière est une fonction gamma incomplète dont nous utiliserons quelques propriétés, rappelées ou démontrées au \S\ref{t202}. Le \S\ref{t137} constitue un intermède méthodologique. Nous y décrivons une seconde démonstration de la relation \eqref{u5}, proposée par l'arbitre anonyme de cet article, et présentant l'intérêt de n'utiliser que la majoration élémentaire de Dirichlet, $\Delta(x)=O(\sqrt{x})$, alors que la démonstration du \S\ref{t135} s'appuie sur une estimation en $O(x^{\delta})$ avec $\delta<1/2$, par exemple l'estimation \eqref{t127} de Voronoï.

\smallskip

La représentation \eqref{u5} de $I(x,v)$ par une intégrale sur la droite critique permet ensuite, au \S\ref{t126}, de mettre en évidence la convergence de $xI(x,v)+\overline{I(1/x,x^2v)}$ vers une fonction continue de $x$ (se prolongeant par continuité en $0$) quand $v$ tend vers l'infini (proposition \ref{t100}), et \eqref{t74} s'en déduit. 

\smallskip

Dans un second temps, nous montrons au \S\ref{t501} que, pour tout $x>0$, la limite quand $v$ tend vers l'infini de la quantité
$$
 \psi(x,v)-x\overline{\psi(1/x,x^2v)}-\demi\log^2\frac 1x-(\gamma-\log 2\pi +\demi \pi i)\log \frac 1x
$$
(limite dont l'existence aura donc été démontrée au \S\ref{t126}) a pour parties réelle et imaginaire les seconds membres de \eqref{t111} et \eqref{t107}. Pour la partie imaginaire, cela découle assez directement des résultats de \cite{baez-duarte-all} : nous en donnons deux démonstrations au \S\ref{t110}. L'étude de la partie réelle requiert en revanche quelques considérations supplémentaires, reposant en particulier sur l'utilisation de la transformation de Hilbert. Nous en rappelons au \S\ref{t115} les propriétés utilisées au \S\ref{t112} pour obtenir \eqref{t111}. 

\smallskip

Notons que nous n'utilisons à aucun moment les fonctions de Bessel, à l'inverse de \cite{Wilton}. 

Nous emploierons à plusieurs reprises la notation d'Iverson : $[P]=1$ si la propriété $P$ est vérifiée, $[P]=0$ sinon.

\section{Rappels sur la fonction $\zeta$ de Riemann}\label{t75}

\subsection{\'Equation fonctionnelle}\label{t82}

Il s'agit de l'identité classique (cf. \cite{titchmarsh-zeta}, chapter II) 
\begin{equation*}
\frac{\zeta(1-s)}{\zeta(s)}=2(2\pi)^{-s}\cos (\pi s/2)\Gamma(s)
\end{equation*}
que nous utiliserons également sous la forme
\begin{equation}\label{t78}
-\pi\frac{\zeta(-s)}{\zeta(s+1)}=(2\pi)^{-s}\sin (\pi s/2)\Gamma(s+1).
\end{equation}
%
Par ailleurs, rappelons la relation $\overline{\zeta(s)}= \zeta(\overline{s})$ pour tout $s\not =1$ qui découle du principe de réflexion de Schwarz. 

\subsection{Majoration dans la bande critique} 
Nous aurons l'usage d'une majoration uniforme pour la fonction $\zeta$ dans la bande critique : pour $0<\sigma < 1$, $\tau\ge 2$, on a
\begin{equation}\label{t77}
 \zeta(\sigma+i\tau) \ll \tau^{(1-\sigma)/2} \log\tau
\end{equation}
 (cf. \cite{ivic}, theorem 1.9 p. 25 par exemple).

\subsection{Moyenne quadratique sur la droite critique}\label{moyquad}

La recherche d'estimations du type  
\begin{equation}\label{t76}
\zeta\left(\demi+i\tau\right)\ll (1+|\tau|)^{\delta} \quad (\tau \in \Real)
\end{equation}
est un des problèmes fondamentaux de la théorie de la fonction $\zeta$ (cf. \cite{titchmarsh-zeta}, chapter V). L'hypothèse de Lindelöf, toujours ouverte, affirme que \eqref{t76} est vraie pour tout $\delta>0$. L'estimation n'est actuellement démontrée que pour $\delta>32/205$ (cf. \cite{huxley-2005} ; dans la prépublication récente \cite{bourgain-2014}, Bourgain démontre \eqref{t76} pour $\delta>53/342$).

Cela étant, l'hypothèse de Lindelöf est {\og vraie en moyenne quadratique\fg} puisque la fonction
$$
I_2(T)=\int_0^T|\zeta(1/2+i\tau)|^2d\tau \quad (T>0)
$$
vérifie l'estimation $I_2(T) \ll T\log (2+T)$.
On connait en fait un développement asymptotique de $I_2(T)$. Posons
$$
I_2(T)=T\log T-(\log 2\pi +1-2\gamma)T+E(T) \quad (T>0).
$$
La fonction $E(T)$ est l'objet d'une abondante littérature (cf. par exemple \cite{ivic}, chapters 4, 15). Nous nous contenterons de l'estimation classique d'Ingham (cf. \cite{53.0313.01}, Theorem A', p. 294) : 
\begin{equation}
\label{eq:ingham}
E(T)\ll T^{1/2}\log (2+T).
\end{equation}  
Dans la proposition suivante, nous utilisons cette estimation pour majorer une intégrale intervenant  au  \S\ref{t126}.

\begin{prop}\label{t88}
Pour $V>0$, on a
$$
\int_{0}^{\infty}\frac{|\zeta(\demi+i\tau)|^2}{1+\tau}\min(1,V^{1/2}\big\vert V-\tau\big\vert^{-1})d\tau \ll V^{-1/2}\log^2 (2+V).
$$
\end{prop}
\dem

Si $V<4$, le minimum figurant dans l'intégrale est $\ll (1+\tau)^{-1}$, donc l'estimation résulte de la 
convergence de l'intégrale
$$
\int_0^{\infty}\frac{|\zeta(\demi+i\tau)|^2}{(1+\tau)^2}\, d\tau
$$
qui se déduit, par exemple, de \eqref{t77}.

On peut donc supposer $V\soe 4$.
La contribution de l'intervalle $|\tau -V|\ioe \sqrt{V}$ à l'intégrale est 
\begin{align*}
 &\ll V^{-1}\int_{V-\sqrt{V}}^{V+\sqrt{V}}|\zeta(1/2+i\tau)|^2d\tau\\
&\ll V^{-1/2}\log V,
\end{align*}
d'après l'estimation \eqref{eq:ingham}.
La contribution des $\tau \soe 2V$ est 
\begin{align*}
 &\ll V^{1/2}\int_{2V}^{\infty}\frac{|\zeta(1/2+i\tau)|^2}{\tau^2}   d\tau\\
&\ll V^{-1/2}\log V,
\end{align*}
en utilisant simplement l'estimation $I_2(T)\ll T\log T$ pour $T\soe 2$ et une intégration par parties. 
De même la contribution des $\tau \ioe V/2$ est 
\begin{align*}
 &\ll V^{-1/2}\int_0^{V/2}\frac{|\zeta(1/2+i\tau)|^2}{1+\tau}  d\tau\\
&\ll V^{-1/2}\log^2 V.
\end{align*}
Maintenant, la contribution de l'intervalle $V+\sqrt{V}<\tau<2V$ est 
\begin{align*}
 &\ll V^{-1/2}\int_{V+\sqrt{V}}^{2V}\frac{|\zeta(1/2+i\tau)|^2}{\tau -V}   d\tau\\
&\ioe V^{-1/2}\sum_{1\ioe k\ioe \sqrt{V}}\frac{1}{k\sqrt{V}}\int_{V+k\sqrt{V}}^{V+(k+1)\sqrt{V}}|\zeta(1/2+i\tau)|^2 d\tau\\
&\ll V^{-1/2}(\log V)\sum_{1\ioe k\ioe \sqrt{V}}\frac 1k\expli{d'après l'estimation \eqref{eq:ingham}}\\
&\ll V^{-1/2}\log^2 V,
\end{align*}
et on a encore la même estimation pour la contribution de l'intervalle $V/2<\tau<V-\sqrt{V}$.\fin


\subsection{Lien avec la fonction arithmétique {\og nombre de diviseurs\fg}}

La fonction arithmétique $\tau(n)$ est liée à la fonction $\zeta$ par la relation
$$
\sum_{n\soe 1}\frac{\tau(n)}{n^s}=\zeta^2(s) \quad (\Re s>1).
$$

D'autre part, on a l'estimation $\tau(n)\ll_{\eps} n^{\eps}$ pour tout $\eps>0$ (cf. \cite{MR1366197}, corollaire 5.3).

\section{Rappels et compléments sur les transformations de Mellin et de Hilbert}

Nous recommandons la lecture de l'appendice A (p. 231) de \cite{baez-duarte-all} et nous en rappelons ci-dessous quelques éléments qui nous seront utiles.

En désignant par $s$ la variable complexe, on note $\sigma=\Re s$ et $\tau=\Im s$. Si $-\infty\ioe a<b\ioe\infty$, on note $\Wcal(a,b)$ l'ensemble des fonctions complexes $f$ mesurables sur $]0,\infty[$ telles que
$$
\int_0^{\infty}|f(t)|t^{\sigma-1}dt<\infty \quad (a<\sigma<b).
$$
Si $f\in \Wcal(a,b)$, la transformée de Mellin $Mf$ définie par \eqref{t123} est holomorphe dans la bande $a<\sigma<b$. 

Nous considérons aux \S\S\ref{t96}, \ref{t201} et \ref{t114} les transformées de Mellin qui interviennent dans notre argumentation, puis nous rappelons d'une part, au \S\ref{t135}, la forme que prend la théorie de Plancherel dans le contexte de la transformation de Mellin, et d'autre part, au \S\ref{t115}, quelques éléments de la théorie de la transformation de Hilbert.

\subsection{La fonction $A$}\label{t96}

Rappelons (cf. \cite{baez-duarte-all}, proposition 10) que $A\in \Wcal(-1,0)$ et que
$$
MA(s)=-\frac{\zeta(-s)\zeta(s+1)}{s(s+1)} \quad (-1<\sigma <0).
$$

On en déduit que la fonction $A_1$ définie par $A_1(t)=A(t)/t=A(1/t)$ appartient à $\Wcal(0,1)$ et que
\begin{equation}\label{t83}
MA_1(s)=\frac{\zeta(1-s)\zeta(s)}{s(1-s)} \quad (0<\sigma <1).
\end{equation}
%
%

\subsection{Le reste dans le problème des diviseurs de Dirichlet}\label{t201}

Nous utiliserons l'expression de la transformée de Mellin du reste $\Delta$, défini par \eqref{t127}, dans le problème des diviseurs de Dirichlet. \`A cette occasion, donnons l'énoncé général d'un principe classique, qui est une réciproque de la proposition 14 de \cite{baez-duarte-all}.
\begin{prop}\label{t129}
Soit $a<b\ioe c<d$, $f\in \Wcal(a,b)$, $g\in\Wcal(c,d)$. On suppose que $g-f$ est un polynôme généralisé
$$
P(t)=\sum_{\rho,k} c_{\rho,k}t^{-\rho}\log^kt,
$$
où la somme est finie, les $\rho$ sont des nombres complexes vérifiant $b\ioe \Re \rho\ioe c$, les $k$ sont des entiers naturels, et les $c_{\rho,k}$ des coefficients complexes.
Alors $Mf$ et $Mg$ sont les restrictions aux bandes $a<\sigma<b$ et $c<\sigma<d$, respectivement, d'une même fonction méromorphe dans la bande $a<\sigma<d$, dont la somme des parties polaires est
$$
\sum_{\rho,k} c_{\rho,k}\frac{(-1)^kk!}{(s-\rho)^{k+1}}.
$$
\end{prop}
\dem

On va montrer que les fonctions $Mf$ et $Mg$ sont des prolongements méromorphes l'une de l'autre.
On a
\begin{equation}\label{t128}
f(t)+P(t)[t\soe 1]=g(t)-P(t)[t< 1].
\end{equation}
En notant $h(t)$ la valeur commune des deux membres de \eqref{t128}, on voit que $h\in \Wcal(a,b)$ (premier membre) et $h\in\Wcal(c,d)$ (second membre). On a donc $h\in \Wcal(a,d)$ ; la transformée de Mellin $Mh$ est holomorphe dans la bande $a<\sigma<d$, et
\begin{align*}
Mf(s) &= Mh(s) +\sum_{\rho,k} c_{\rho,k}\frac{(-1)^kk!}{(s-\rho)^{k+1}} \quad (a<\sigma<b),\\
Mg(s) &= Mh(s) +\sum_{\rho,k} c_{\rho,k}\frac{(-1)^kk!}{(s-\rho)^{k+1}} \quad (c<\sigma<d).\fine
\end{align*}

La proposition \ref{t129} s'applique aux fonctions 
\begin{align*}
f(t) =\sum_{n\ioe t}\tau(n), \, 
g(t) =\Delta(t) \textrm{ et } 
P(t)=-t(\log t +2\gamma-1).
\end{align*}
En utilisant l'estimation $\Delta(t)=O(t^{\delta})$ (où $\delta=\demi$ si l'on se contente de l'estimation élémentaire de Dirichlet, et $\delta=\frac 13$ d'après l'estimation de Voronoï \eqref{t127}), on a $f\in\Wcal(-\infty,-1)$, $g\in\Wcal(-1,-\delta)$. 
Comme
\begin{equation*}
Mf(s)=\int_0^{\infty}\big (\sum_{n\ioe t}\tau(n)\big)t^{s-1}dt
=\sum_{n\soe 1}\tau(n)\int_n^{\infty}t^{s-1}dt
=\frac{\zeta^2(-s)}{-s}\quad(\sigma<-1),
\end{equation*}
on en déduit que 
$$
M\Delta(s)=\frac{\zeta^2(-s)}{-s}\quad(-1<\sigma<-\delta).
$$
Par conséquent la fonction $\Delta_1$ définie par $\Delta_1(t)=\Delta(t)/t$ appartient à $\Wcal(0,1-\delta)$ et 
\begin{equation}\label{t84}
M\Delta_1(s)=\frac{\zeta^2(1-s)}{1-s}\quad(0<\sigma<1-\delta).
\end{equation}

\subsection{La fonction $t\mapsto e^{it}[t\ioe v]$}\label{t114}

Pour tout $v>0$, la fonction $t\mapsto e^{it}[t\ioe v]$ appartient à $\Wcal(0,\infty)$ et sa transformée de Mellin est $\lambda(s,v)$ d\'efinie par \eqref{t99}. Par conséquent, pour $x>0$ et $v>0$, la transformée de Mellin de $t\mapsto e^{2\pi it x}[t\ioe v]$ est $(2\pi x)^{-s}\lambda(s,2\pi xv)$.

\subsection{La transformation de Mellin-Plancherel sur $L^2$}\label{t135}

Si $f\in L^2(0,\infty)=L^2$, la formule \eqref{t123}, où l'intégrale doit être comprise comme $\lim_{T\vers \infty}\int_{1/T}^T$ dans $\Lcal^2=L^2(\demi+i\Real,d\tau/2\pi)$, définit un élément de ce dernier espace, et l'application ainsi définie, dite transformation de Mellin-Plancherel, est une isométrie bijective entre espaces de Hilbert (cf. par exemple \cite{titchmarsh-fourier} \S 3.17 pp. 94-95). En particulier, si $f$ et $g$ appartiennent à $L^2$, on a
$$
\int_0^{\infty}f(t)g(t)\,dt =\int_{\sigma=1/2}Mf(s)Mg(1-s)\,\frac{d\tau}{2\pi}
$$
(théorème de Plancherel), où les deux intégrales sont absolument convergentes.

Appliquons maintenant le théorème de Plancherel à l'intégrale $I(x,v)$ définie par \eqref{u136}. En utilisant le fait que $\Delta_1 \in L^2$, on obtient
\begin{align*}
I(x,v)&=2\int_0^{\infty}\Delta_1(t)\cdot e^{2\pi itx}[t\ioe v]dt\\
&=\int_{\sigma=1/2}\frac{\zeta^2(1-s)}{1-s}\lambda(1-s,2\pi xv)(2\pi x)^{s-1}\,\frac{d\tau}{\pi}\\
&=\int_{\sigma=1/2}\frac{\zeta^2(s)}{s}\lambda(s,2\pi xv)(2\pi x)^{-s}\,\frac{d\tau}{\pi},
\end{align*}
ce qui démontre la relation \eqref{u5}.

\subsection{La transformation de Hilbert}\label{t115}

Nous donnons dans ce paragraphe les faits que nous utiliserons concernant cette transformation, objet du chapitre V de \cite{titchmarsh-fourier}, ainsi que du traité très complet \cite{king-09}. La transformée de Hilbert d'une fonction $f$ définie presque partout sur $]-\infty,\infty[$ et mesurable est
$$
Hf(x)=\frac 1{\pi}\vp\int_{-\infty}^{\infty}\frac{f(t)}{x-t}dt
$$
où $\vp$ signifie {\og valeur principale\fg} (de Cauchy) :
\begin{equation}\label{t116}
Hf(x)=\frac 1{\pi}\lim_{\delta\vers 0}\int_{\vert x-t\vert >\delta}\frac{f(t)}{x-t}dt.
\end{equation}


Les propriétés fondamentales de la transformation de Hilbert sont, d'une part, le couple d'équations
\begin{equation*}
H\cos =\sin \quad ; \quad H\sin=-\cos.
\end{equation*}
et, d'autre part, le fait que $H$ commute avec les opérateurs transformant la fonction $f(x)$ en $f(\lambda x)$ et $f(x+c)$ ($\lambda >0$, $c\in \Real$). 

De plus, si l'on étend la définition \eqref{t116} en remplaçant la condition $\vert x-t\vert >\delta$ par $1/\delta >\vert x-t\vert >\delta$, la transformée de Hilbert d'une constante est nulle. Par conséquent, $H$ transforme formellement une série trigonométrique en la série trigonométrique \emph{conjuguée} (au sens de \cite{MR1963498}, chapter I, ($1\cdot 3$)). 

\smallskip

Si $f\in L^2(\Real)$, la limite \eqref{t116} existe au sens de la topologie de cet espace de Hilbert, et l'opérateur ainsi défini est unitaire.  

Dans l'espace $L^2(0,\infty)$ (que nous avons noté simplement $L^2$), la transformation de Hilbert a deux versions, paire et impaire :
\begin{align}
H_0f(x) &=\frac{1}{\pi}\vp\int_0^{\infty}f(t) \Big (\frac{1}{x-t}+\frac{1}{x+t}\Big)dt\label{t130}\\ 
H_1f(x) &=\frac{1}{\pi}\vp\int_0^{\infty}f(t) \Big (\frac{1}{x-t}-\frac{1}{x+t}\Big)dt \label{t132}
\end{align}
Ce sont deux opérateurs unitaires sur $L^2$, commutant avec les dilatations, et donc \emph{diagonalisables} via la transformation de Mellin-Plancherel. Pour $f\in L^2$, on a les identités suivantes\footnote{Les calculs menant aux équations (5.126) et (5.129), p. 274 de \cite{king-09}, sont valables si $f$ est indéfiniment dérivable et à support compact dans $]0,\infty[$, et les deux identités \eqref{t124} et \eqref{t117} s'en déduisent par densité.} entre éléments de $\Lcal^2$ :
\begin{align}
MH_0f(s)&=\tan (\pi s/2)Mf(s)\label{t124}\\
MH_1f(s)&=-\cot (\pi s/2)Mf(s),\label{t117}
\end{align}
Les fonctions $\tan (\pi s/2)$ et $\cot (\pi s/2)$ sont de module $1$ sur la droite $\sigma=1/2$ ; nous utiliserons également le fait qu'elles sont bornées dans toute bande verticale fermée incluse dans la bande $0<\sigma<1$.


\smallskip

Les conditions suivantes sont suffisantes pour que \eqref{t130} et \eqref{t132} soient bien définies au point $x_0>0$ :

$\bullet$ la fonction $t\mapsto f(t)/(1+\vert t\vert)$ est intégrable sur $]0,\infty[$ ;

$\bullet$ il existe $a$ tel que la fonction $t\mapsto \big(f(x_0-t)-a\big)/t$ soit intégrable au voisinage de $0$ (condition de Dini).

De plus, si $f$ est continue au voisinage de $x_0$ et s'il existe $h_1,h_2>0$ et $\fhi\in L^1(-h_2,h_2)$ tels que
$$
\vert \big(f(x-t)-f(x) \big)/t\vert\ioe \fhi(t) \quad (\vert x-x_0\vert \ioe h_1,\; \vert t\vert \ioe h_2),
$$
alors $H_0$ et $H_1$ sont continues en $x_0$.

\subsection{La fonction $B$}

Dans cet article, nous utiliserons la transformée de Hilbert impaire de la fonction $A_1$ (multipliée par $\pi$) :
$$
B(x)=\vp\int_0^{\infty}A(t) \Big (\frac{1}{x-t}-\frac{1}{x+t}\Big)\, \frac{dt}{t} \quad (x>0)
$$
Notons que $A_1(t)/(1+\vert t\vert)$ est intégrable et que, pour tout $x_0>0$, il existe $h_1$ et $h_2$ positifs tels que 
$$
A_1(x-t)-A_1(x)\ll |t|\log(1/|t|) \quad (|x-x_0|<h_1, \; 0<|t|<h_2),
$$
(cf. \cite{baez-duarte-all}, propositions 1 et 7). On en déduit que $B(x)$ existe pour tout $x>0$ et est une fonction continue de $x$. La proposition suivante donne quelques identités vérifiées par $B$ ; nous omettons les démonstrations, qui sont de simples manipulations algébriques (tenant compte notamment de la relation $A(x)=xA(1/x)$).

\begin{prop}
Pour $x>0$, on a
\begin{align}
B(1/x)&=x\cdot\vp\int_0^{\infty}A(t) \Big (\frac{-1}{x-t}-\frac{1}{x+t}\Big)\, \frac{dt}{t}\notag\\
xB(x)+B(1/x)&=-2x\int_0^{\infty}\frac{A(t)}{x+t}\, \frac{dt}{t}\label{t133}\\
\frac{d}{dx}\big(xB(x)+B(1/x)\big) &=-2\int_0^{\infty}\frac{A(t)}{(x+t)^2}\, dt\\
B(1/x)-x\frac{d}{dx}\big(xB(x)+B(1/x)\big)&=-2x\cdot\vp\int_0^{\infty}A(t)\frac{x^2+t^2}{(x-t)(t+x)^2}\, \frac{dt}{t}.\label{t134}
\end{align}
\end{prop}

\smallskip

Comme $A_1\in L^2$, on a également $B\in L^2$ et les identités \eqref{t83} et \eqref{t117} fournissent l'expression suivante de la transformée de Mellin-Plancherel de $B$ : on a pour presque tout $s$ tel que $\sigma=1/2$,
\begin{equation}\label{t125}
MB(s)=-\pi\cot(\pi s/2)\cdot\frac{\zeta(s)\zeta(1-s)}{s(1-s)}.
\end{equation}
En notant $F(s)$ la fonction méromorphe sur $\Com$ figurant au second membre de \eqref{t125}, on constate que \eqref{t77} entraîne l'estimation $F(s)\ll \log^2(2+\lvert s\rvert)(1+\vert s\vert)^{-3/2}$ uniformément dans toute bande verticale fermée incluse dans la bande $0<\sigma<1$. On peut donc écrire par inversion de Mellin
\begin{equation}\label{t119}
B(x)=\frac{1}{2\pi i}\int_{\sigma=a}F(s)x^{-s} ds,
\end{equation}
où $a$ est choisi arbitrairement tel que $0<a<1$. La relation \eqref{t119}, a priori valable pour presque tout $x>0$, l'est sans exception, par continuité. Sur cette expression, on voit que $B$ est uniformément $O_a(x^{-a})$ sur $]0,\infty[$ pour chaque $a$ tel que $0<a<1$, d'où il découle que $B\in \Wcal(0,1)$. On en déduit que
\begin{equation}\label{t120}
MB(s)=F(s) \quad (0<\sigma<1).
\end{equation}

\section{Résultats auxiliaires sur les fonction gamma incomplètes et les fonctions cosinus et sinus intégral généralisées}\label{t202}

Dans ce paragraphe, nous considérons les fonctions
\begin{align*}
\ci(s,v)=\int_v^{\infty}\cos t \cdot t^{s-1}dt\, ; &\quad \Ci(s,v)=\int_0^{v}\cos t \cdot t^{s-1}dt\\
\si(s,v)=\int_v^{\infty}\sin t \cdot t^{s-1}dt\, ; &\quad \Si(s,v)=\int_0^{v}\sin t \cdot t^{s-1}dt\\
\Lambda(s,v)=\int_v^{\infty}e^{it} \cdot t^{s-1}dt\, ; &\quad \lambda(s,v)=\int_0^{v}e^{it}\cdot t^{s-1}dt
\end{align*}
pour $v>0$ et $0<\Re s <1$.
Notons pour commencer que $\Ci$, $\Si$ et $\lambda$ sont absolument convergentes et que $\ci$, $\si$ et $\Lambda$ sont semi-convergentes. D'autre part, les relations classiques
\begin{equation} 
\label{eq:mellin-trigo}
\int_0^{\infty}\cos t \cdot t^{s-1}dt\\
=\Gamma(s)\cos(\pi s/2) \,\textrm{ et}\,  \int_0^{\infty}\sin t \cdot t^{s-1}dt\\
=\Gamma(s)\sin(\pi s/2) \quad(0<\Re s <1)
\end{equation}
(cf. \cite{37.0450.01}, \S 62, (4) et (5)) et l'équation fonctionnelle de la fonction $\zeta$ entraînent 
\begin{align}
2(2\pi)^{-s}\big(\ci(s,v)+\Ci(s,v)\big)&=2(2\pi)^{-s}\Gamma(s)\cos(\pi s/2)\notag\\
&=\frac{\zeta(1-s)}{\zeta(s)}\label{t143}
\end{align}
et
\begin{align}
2(2\pi)^{-s}\big(\si(s,v)+\Si(s,v)\big)&=2(2\pi)^{-s}\Gamma(s)\sin(\pi s/2)\notag\\
&=\tan(\pi s/2)\frac{\zeta(1-s)}{\zeta(s)}.\label{u143}
\end{align}
Par conséquent,
\begin{align}
2(2\pi)^{-s}\big(\lambda(s,v)+\Lambda(s,v)\big)&=2(2\pi)^{-s}\Gamma(s)e^{i\pi s/2}\notag\\
&=G(s)\frac{\zeta(1-s)}{\zeta(s)},\label{v143}
\end{align}
avec 
\[
G(s)= 1+i \tan(\pi s/2). 
\] 
On déduit de \eqref{v143}  l'identité suivante, qui interviendra dans la démonstration de la proposition~\ref{t100}.
\begin{prop}\label{t101}
Pour $x>0$, $v>0$ et $s\in \Com$ tel que $\Re s =\demi$, on a
\begin{align*}
\frac{\zeta^2(s)}{s}(2\pi )^{-s}\big(x^{1-s}\lambda(s,v)+x^{s}\overline{\lambda(1-s,v)}\,\big)
&=
\frac{\zeta(s)\zeta(1-s)}{2s}(x^{1-s}G(s)+x^s\overline{G(1-s)})
\\ &\quad - \frac{\zeta^2(s)}{s}(2\pi )^{-s}\big(x^{1-s}\Lambda(s,v)+x^{s}\overline{\Lambda(1-s,v)}\,\big).
\end{align*}
\end{prop}

\medskip

Les estimations que nous utiliserons seront démontrées grâce à la proposition classique suivante, qui découle d'une intégration par parties en écrivant $ge^{if}=(g/f')\cdot f'e^{if}$.

\begin{prop}\label{t12}
Soit $a$ et $b$ deux nombres réels, $a<b$, et soit $f,g \, : [a,b] \vers \Real$ deux fonctions continûment dérivables telles que 

$\bullet$ la fonction $f'$ ne s'annule pas ;

$\bullet$ la fonction $g/f'$ est monotone et de signe constant ;

$\bullet$ on a $\lvert g(t)/f'(t)\rvert \ioe c$ pour $a<t<b$.

Alors
$$
\Big\lvert\int_a^b g(t)e^{if(t)}\, dt \Big\rvert \ioe 2c.
$$
\end{prop}

La proposition suivante rassemble les estimations dont nous aurons l'usage pour les fonctions $\lambda$ et $\Lambda$ dans la preuve de la proposition \ref{prop:eq-approchee-psi1}. 
\begin{prop}\label{t131}
Pour $v>0$ et $\tau \in \Real$, on a
\begin{align*}
|\lambda(1/2+i\tau,v)|&\ioe \min\big(4,2\sqrt{v}/(|\tau|-v)\big) \quad (|\tau|>v)\\
|\Lambda(1/2+i\tau,v)|&\ioe \min\big(4,2\sqrt{v}/(v- |\tau|)\big) \quad (|\tau|<v).
\end{align*}
\end{prop}
\dem

Nous donnons la démonstration pour $\lambda$, celle pour $\Lambda$ étant similaire. 

On suppose donc  $\lvert\tau \vert>v$. On a
\begin{equation*}\label{t17}
\lambda(1/2+i\tau,v)=\int_0^v\frac{e^{i(t+\tau\log t)}}{\sqrt{t}}dt 
\end{equation*}

Avec les notations de la proposition \ref{t12}, on a ici
\begin{align*}
g(t) &=t^{-1/2},\\
f(t) &=t+\tau\log t,
\end{align*}
de sorte que
$$
f'(t)/g(t)= \sqrt{t}+\tau/\sqrt{t}.
$$
Cette fonction est monotone sur $]0,\lvert\tau\rvert]$, donc sur $]0,v]$, et on a $|f'(t)/g(t)|\soe (\lvert\tau\rvert-v)/\sqrt{v}$ sur ce dernier intervalle. Par conséquent,
\begin{equation}\label{t15}
\Big \vert\int_0^v\frac{e^{i(t+\tau\log t)}}{\sqrt{t}}dt\Big \vert \ioe \frac{2\sqrt{v}}{\lvert\tau\rvert-v}.
\end{equation}
Si $v$ et $\lvert\tau\rvert$ sont trop proches, on peut améliorer cette inégalité de la façon suivante. Si $\lvert\tau\rvert \soe v+\sqrt{v}$, la majoration \eqref{t15}  donne
\begin{equation}\label{t86}
|\lambda(1/2+i\tau,v)| \ioe 2.
\end{equation}
En revanche, si $\lvert\tau\rvert-\sqrt{v}<v<\lvert\tau\rvert$, deux sous-cas se présentent :

$\bullet$ si $\lvert\tau\rvert-\sqrt{v} < 0$, alors $v<\lvert\tau\rvert<\sqrt{v}$, donc $v<1$ et
$$
|\lambda(1/2+i\tau,v)| \ioe \int_{0}^1\frac{dt}{\sqrt{t}}=2 \, ;
$$

$\bullet$ si $\lvert\tau\rvert-\sqrt{v}>0$, alors
$$
|\lambda(1/2+i\tau,v)| \ioe \big|\lambda(1/2+i\tau,|\tau|-\sqrt{v})\big| +\int_{|\tau|-\sqrt{v}}^v\frac{dt}{\sqrt{t}}.
$$
Comme
$$
\lvert\tau\rvert \soe \lvert\tau\rvert-\sqrt{v}+\sqrt{\lvert\tau\rvert-\sqrt{v}},
$$
on a $|\lambda(1/2+i\tau,\lvert\tau\rvert-\sqrt{v})|\ioe 2$ d'après \eqref{t86}. Enfin,
\begin{align*}
\int_{\lvert\tau\rvert-\sqrt{v}}^v\frac{dt}{\sqrt{t}}&= 2\frac{\sqrt{v}+v-\lvert\tau\rvert}{\sqrt{v}+\sqrt{\lvert\tau\rvert-\sqrt{v}}}\\
&\ioe 2.
\end{align*}
On en déduit la majoration $|\lambda(1/2+i\tau,v)|\ioe 4$.\fin

\section{Seconde démonstration de l'identité \eqref{u5}}\label{t137}

La démonstration de l'identité \eqref{u5} proposée par l'arbitre anonyme part de la constatation suivante : l'estimation de Dirichlet $\Delta(x)=O(\sqrt{x})$ entraîne que, pour tout $\eta$ tel que $0<\eta<1/2$, la fonction
$$
x\mapsto \Delta_1(x)x^{-\eta}=\Delta(x)x^{-1-\eta}
$$
appartient à $L^2(0,\infty)$. Sa transformée de Mellin étant d'après \eqref{t84}
$$
\int_0^{\infty}\Delta_1(t)t^{s-\eta -1} dt=\frac{\zeta^2(1-s+\eta)}{1-s+\eta} \quad (\eta <\sigma <1/2+\eta),
$$
le théorème de Plancherel donne
\begin{align}
I(x,v,\eta)&=2\int_0^{\infty}\Delta_1(t)t^{-\eta}\cdot e^{2\pi itx}[t\ioe v]dt\notag\\
&=\int_{\sigma=1/2}\frac{\zeta^2(1-s+\eta)}{1-s+\eta}\lambda(1-s,2\pi xv)(2\pi x)^{s-1}\,\frac{d\tau}{\pi}\notag\\
&=\int_{\sigma=1/2}\frac{\zeta^2(s+\eta)}{s+\eta}\lambda(s,2\pi xv)(2\pi x)^{-s}\,\frac{d\tau}{\pi}.\label{t138}
\end{align}
Or l'estimation \eqref{t77} et la proposition \ref{t131} fournissent pour tout $v$ fixé l'estimation
$$
\frac{\zeta^2(s+\eta)}{s+\eta}\lambda(s,2\pi xv) \ll \frac{\log^2\tau}{\tau^{3/2}} \quad (\sigma=\demi,\, \lvert \tau \rvert \soe 2, \, 0<\eta<1/2).
$$
Compte tenu du théorème de convergence dominée, on obtient donc \eqref{u5} en faisant tendre $\eta$ vers $0$ dans \eqref{t138}.

\section{Démonstration de la proposition \ref{prop:eq-approchee-psi1}}\label{t126}

Soit $\Fcal$ l'ensemble des fonctions continues $f :\, ]0,\infty[ \vers \Com$ telles que $f(x)$ et $\widetilde{f}(x)=x\overline{f(1/x)}$ se prolongent par continuité en $0$. Les constantes appartiennent à $\Fcal$, ainsi que les fonctions continues sur $]0,\infty[$ vérifiant
\begin{align*}
f(x)&=o(1)\quad (x\vers 0)\\
f(x)&=o(x)\quad (x\vers \infty).
\end{align*}
Observons que $\Fcal$ est un espace vectoriel sur $\Com$, invariant par l'involution $f\mapsto \widetilde{f}$. Dans ce qui suit nous écrirons des égalités du type
\begin{equation}\label{t136}
f_1(x,v)=f_2(x,v)\pmod \Fcal
\end{equation}
pour signifier que la différence $f_1-f_2$ est une fonction de la seule variable $x$, appartenant de plus à $\Fcal$. Dans ce qui suit, il sera crucial de s'assurer que $f_1(x,v)-f_2(x,v)$ ne dépend pas de la variable $v$ afin d'écrire toute relation du type \eqref{t136}.

\medskip

Nous commençons par un calcul d'intégrale.

\begin{prop}\label{t98}
On a
$$
a_0=\int_0^1\frac{e^{it} -1}{t}dt +\int_{1}^{\infty} \frac{e^{it}}{t}dt =-\gamma+i\pi/2.
$$
\end{prop}
\dem

On a, d'après \eqref{eq:mellin-trigo}, pour $0<\Re(s)<1$,
\begin{align*}
e^{i\pi s/2}\Gamma(s) &=\int_0^{\infty}e^{it}t^{s-1}dt=\int_0^1 +\int_1^{\infty}\\
&=\frac 1s+\int_0^1 \big (e^{it}-1\big)t^{s-1}dt+\int_1^{\infty}e^{it}t^{s-1}dt\\
&=\frac 1s+a_0+o(1),
\end{align*}
quand $s$ tend vers $0$ par valeurs positives, et d'autre part,
\begin{align*}
e^{i\pi s/2}\Gamma(s) &=\big(1+i\pi s/2 +O(s^2)\big)\frac{\Gamma(1+s)}{s}\\
&=\frac 1s-\gamma+i\pi/2+o(1),
\end{align*}
puisque $\Gamma'(1)=-\gamma$.\fin

\smallskip

Nous passons maintenant à l'étude d'une intégrale dépendant de $x$ et de $v$.

\begin{prop}\label{t91}
Pour $x>0$ et $v>0$ on a
$$
\int_{1}^v\frac{e^{2\pi it x}}{t}(\log t +2\gamma) dt=\demi\log^2\frac 1x+\Big(\gamma-\log 2\pi +i\frac{\pi}{2}\Big)\log\frac 1x +\eps_0(x,v) \pmod \Fcal
$$
où
$$
\eps_0(x,v)=-\int_v^{\infty}\frac{e^{2\pi it x}}{\pi t}(\log t +2\gamma) dt.
$$
\end{prop}
\dem

On a
$$
\int_1^v\frac{e^{2\pi it x}}{\pi t}(\log t +2\gamma) dt=\int_1^{\infty}-\int_v^{\infty},
$$
où les intégrales sont semi-convergentes. Le dernier terme est $\eps_0(x,v)$.

Maintenant,
\begin{align*}
\int_1^{\infty}\frac{e^{2\pi it x}}{t}(\log t +2\gamma) dt&=\int_{2\pi x}^{\infty} (\log (t/2\pi x) +2\gamma) \frac{e^{it}}{t} dt\\
&=\int_{2\pi x}^{\infty}  \log t \cdot\frac{e^{it}}{t}dt-(\log x+\log 2\pi -2\gamma)\int_{2\pi x}^{\infty} \frac{e^{it}}{t}dt\, ,
\end{align*}
avec
\begin{align*}
\int_{2\pi x}^{\infty}  \log t \cdot\frac{e^{it}}{t}dt &=\int_{2\pi x}^1\log t \frac{dt}{t}+\int_{2\pi x}^1\log t\cdot\frac{e^{it} -1}{t}dt +\int_{1}^{\infty} \log t \cdot\frac{e^{it}}{t}dt\\
&=-\demi\log^2 2\pi x \pmod \Fcal
\end{align*}
et
\begin{align*}
\int_{2\pi x}^{\infty} \frac{e^{it}}{t}dt &=\int_{2\pi x}^1\ \frac{dt}{t}-\int_0^{2\pi x}\frac{e^{it} -1}{t}dt  +\int_0^1\frac{e^{it} -1}{t}dt +\int_{1}^{\infty} \frac{e^{it}}{t}dt\\
&=\log(1/2\pi x) +a_0 -\int_0^{2\pi x}\frac{e^{it} -1}{t}dt,
\end{align*}
où $a_0=-\gamma +i\pi/2$, d'après la proposition \ref{t98}. Comme la fonction
$$
x\mapsto (\log x+\log 2\pi -2\gamma)\int_0^{2\pi x}\frac{e^{it} -1}{t}dt
$$
appartient à $\Fcal$, on en déduit le résultat annoncé.\fin

\smallskip

Rappelons la notation 
$$
\psi(x,v)=\sum_{n\le  v} \frac{\tau(n)}{n} e^{2\pi inx}. 
$$
Dans la proposition suivante, nous ramenons l'étude de $\psi(x,v)$ à celle de l'intégrale $I(x,v)$ définie par \eqref{u136}.

\begin{prop}\label{u140}
Pour $x>0$, $v>0$, on a
\begin{equation}\label{u142}
\psi(x,v)=-i\pi xI(x,v)+\frac{1}{2}\log^2\frac 1x +\Big(\gamma -\log 2\pi +i\frac{\pi}{2}\Big)\log \frac 1x +\eps_1(x,v)\pmod \Fcal,
\end{equation}
où $I(x,v)$ est défini par \eqref{u136}, et
$$
\eps_1(x,v)=-\int_v^{\infty}\frac{e^{2\pi it x}}{\pi t}(\log t +2\gamma) dt+\frac{e^{2\pi ivx}}{v}\Delta(v)-\int_v^{\infty}\Delta(t)\frac{e^{2\pi it x}}{t^2}dt.
$$
\end{prop}
\dem

Nous commençons par effectuer une intégration par parties sur l'expression de $\psi(x,v)$ comme intégrale de Stieltjes :
\begin{align}
\psi(x,v)&=\sum_{n\le  v} \frac{\tau(n)}{n} e^{2\pi in x}
=\int_{1^-}^v\frac{e^{2\pi it x}}{t}d\Big (\sum_{n\ioe t}\tau(n)\Big )\notag\\
&=\int_{1}^v\frac{e^{2\pi it x}}{t}(\log t +2\gamma) dt+\int_{1^-}^v\frac{e^{2\pi it x}}{t}d\Delta(t).\label{u6}
\end{align}
La première intégrale de \eqref{u6} est l'objet de la proposition \ref{t91}. Quant à la seconde, on a :
\begin{align}
\int_{1^-}^v&\frac{e^{2\pi it x}}{t}d\Delta(t)\\&=\frac{e^{2\pi it x}}{t}\Delta(t) \Big\vert_{1^-}^v-\int_1^v\Delta(t)\Big (-\frac{e^{2\pi it x}}{t^2}+2 i \pi x\frac{e^{2\pi it x}}{t}\Big)dt\notag\\
&=(2\gamma -1)e^{2i\pi x}+\int_1^{\infty}\Delta(t)\frac{e^{2\pi it x}}{t^2}dt +2i\pi x\int_0^1\Delta(t)\frac{e^{2\pi it x}}{t}dt-i\pi xI(x,v)+\eps_2(x,v),\label{u9}
\end{align}
où $I(x,v)$ est définie par \eqref{u136} et où l'on a posé
\begin{equation*}
\eps_2(x,v)=\frac{e^{2\pi ivx}}{v}\Delta(v)-\int_v^{\infty}\Delta(t)\frac{e^{2\pi it x}}{t^2}dt.
\end{equation*}
Les trois premiers termes de \eqref{u9} sont des fonctions de la seule variable $x$, appartenant à $\Fcal$ (pour le troisième, cela résulte du lemme de Riemann-Lebesgue). 
On en déduit le résultat annoncé, avec $\eps_1=\eps_0+\eps_2$.\fin

\medskip

La proposition suivante contient l'argument principal de la démonstration, à savoir l'approximation de $xI(x,v)+\overline{I(1/x,x^2v)}$ par une transformée de Mellin inverse absolument convergente.

\begin{prop}\label{t100}
Pour $x>0$, $v>0$, on a
\begin{equation}\label{u144}
xI(x,v)+\overline{I(1/x,x^2v)}=\eta(x,v) \pmod \Fcal\, ,
\end{equation}
où
\begin{equation}\label{u147}
\begin{aligned}
\eta(x,v)&=-\int_{|\tau|\ioe V}\frac{\zeta^2(s)}{s}(2\pi )^{-s}\big(x^{1-s}\Lambda(s,V)+x^{s}\overline{\Lambda(1-s,V)}\,\big)\frac{d\tau}{\pi}\\
&\qquad+ 
\int_{|\tau|> V}\frac{\zeta^2(s)}{s}(2\pi )^{-s}\big(x^{1-s}\lambda(s,V)+x^{s}\overline{\lambda(1-s,V)}\,\big)\,\,\frac{d\tau}{\pi}\\
&\qquad\qquad-\int_{|\tau|>V}\zeta(s)\zeta(1-s)H(s) x^s\frac{d\tau}{2\pi},
\end{aligned}
\end{equation}
avec $V=2\pi xv$, $s=\demi +i\tau$ 
et $\displaystyle{H(s)=\frac {1}{s(1-s)}+i\left(\frac{\cot\pi s/2}{1-s}-\frac{\tan\pi s/2}{s}\right).}$ 
\end{prop}
\dem

Appliquons  l'identité \eqref{u5} aux couples $(x,v)$ et $(1/x,x^2v)$ :
\begin{align*}
I(x,v) &=\int_{\Re s=\demi}\frac{\zeta^2(s)}{s}\lambda(s,2\pi xv)(2\pi x)^{-s}\,\frac{d\tau}{\pi}\\
I(1/x,x^2v)&=\int_{\Re s=\demi}\frac{\zeta^2(1-s)}{1-s}\lambda(1-s,2\pi xv)(2\pi /x)^{s-1}\,\frac{d\tau}{\pi},
\end{align*}
donc
\begin{align*}
xI(x,v)+\overline{I(1/x,x^2v)}&= \int_{\Re s=\demi}\frac{\zeta^2(s)}{s}(2\pi )^{-s}\big(x^{1-s}\lambda(s,V)+x^{s}\overline{\lambda(1-s,V)}\,\big)\,\frac{d\tau}{\pi}\\
&=\int_{|\tau|\ioe V}+\int_{|\tau|> V}\, ,
\end{align*}
où $V=2\pi xv$.

L'intégrale $\int_{|\tau|> V}$ est la deuxième intervenant dans $\eta(x,v)$. Quant à l'intégrale $\int_{|\tau|\ioe V}$, elle vaut d'après la proposition \ref{t101}

\begin{align}
\int_{|\tau|\ioe V}&\frac{\zeta(s)\zeta(1-s)}{s}(x^{1-s}G(s)+x^s\overline{G(1-s)})\,\frac{d\tau}{2\pi}\label{t102}
\\&-\int_{|\tau|\ioe V}\frac{\zeta^2(s)}{s}(2\pi )^{-s}\big(x^{1-s}\Lambda(s,V)+x^{s}\overline{\Lambda(1-s,V)}\,\big)\frac{d\tau}{\pi}.\label{t102bis}
\end{align}
L'intégrale \eqref{t102bis} est la première intervenant dans $\eta(x,v)$. Nous réorganisons l'intégrale \eqref{t102} : en séparant les termes correspondant aux expressions $x^s$ et $x^{1-s}$, en effectuant un changement de variable $s\mapsto 1-s$ dans ces derniers, et finalement en regroupant, on obtient l'intégrale
$$
\int_{|\tau|\ioe V}\zeta(s)\zeta(1-s)H(s)x^s\,\frac{d\tau}{2\pi},
$$
où
\begin{equation*}
\label{defH}
H(s) =\frac{G(1-s)}{1-s}+\frac{\overline{G(1-s)}}{s}= \frac 1{s(1-s)}+i\Big (\frac{\cot\pi s/2}{1-s}-\frac{\tan\pi s/2}{s}\Big)
\end{equation*}
puisque $\overline{\tan z}=\tan{\overline{z}}$ et $\overline{s}=1-s$ si $\Re s =\demi$.
Par ailleurs l'identité 
\begin{equation*} 
\label{eq:th}
\tan\pi s/2=\frac{1+i \tanh(\pi \tau/2)}{1-i \tanh(\pi \tau/2)} \quad(s=1/2+i\tau, \tau\in\Real).
\end{equation*}
montre que $\tan(\pi s/2)= \pm i + o(\tau^{-1}) \,(s= \frac12 +i\tau, \tau\to \pm\infty)$, et cela fournit la majoration  
\begin{equation}\label{eq:majH} 
H(s) \ll (1+\tau^2)^{-1} \quad(s= \frac12 + i\tau). 
\end{equation}
La fonction $s\mapsto\zeta(s)\zeta(1-s)H(s)x^s$ est donc intégrable sur la droite $\Re s=1/2$, et nous pouvons ainsi écrire
\begin{equation}\label{t104}
\int_{|\tau|\ioe V}\zeta(s)\zeta(1-s)H(s)x^s\,\frac{d\tau}{2\pi}=\int_{-\infty}^{\infty}-\int_{|\tau|>V}.
\end{equation}
La dernière intégrale de \eqref{t104} est la troisième composant $\eta(x,v)$. Quant à l'avant-dernière, c'est une fonction continue de $x>0$ qui est uniformément $O(x^{1/2})$ et appartient donc à $\Fcal$, d'où le résultat.\fin

\medskip

Pour conclure la démonstration de la proposition \ref{prop:eq-approchee-psi1}, il nous reste à estimer les fonctions $\eps_1$ et $\eta$. C'est l'objet des deux propositions suivantes.
\begin{prop}\label{u146}
Soit $K_1,K_2>0$. Pour $0<x\ioe K_1$, $v>0$ et $x^2v\soe K_2$, on a
$$
\eps_1(x,v)-x \overline{\eps_1(1/x,x^2v)} \ll (x^2v)^{-1/2},
$$
où la constante implicite dépend uniquement de $K_1$ et $K_2$.
\end{prop}
\dem

Nous commençons par estimer $\eps_1(x,v)$ sous la seule hypothèse $v\soe K$ avec $K>0$. Rappelons que
$$
\eps_1(x,v)=-\int_v^{\infty}\frac{e^{2\pi it x}}{\pi t}(\log t +2\gamma) dt+\frac{e^{2\pi ivx}}{v}\Delta(v)-\int_v^{\infty}\Delta(t)\frac{e^{2\pi it x}}{t^2}dt.
$$
Par le second théorème de la moyenne, la première intégrale est $\ll \big(\log (2+v)\big)/xv$. En utilisant l'estimation de Dirichlet $\Delta(v) \ll v^{1/2}$, on voit que les autres termes sont $\ll v^{-1/2}$.
On en déduit, si $0<x\ioe K_1$ et $x^2v\soe K_2$ :
\begin{align*}
\eps_1(x,v)-x\overline{\eps_1(1/x,x^2v)} &\ll \big(\log (2+v)\big)/xv +v^{-1/2} +x\big(\log(2+x^2v)\big)/xv +x(x^2v)^{-1/2}\\
&\ll (x^2v)^{-1/2}.\fine
\end{align*}

\begin{prop}\label{prop:eta} 
Pour $0<x\ioe K_1$, $v>0$ et $x^2v\soe K_2$, on a
$$
\eta(x,v) \ll (x^2v)^{-1/2}\log^2(2+x^2v),
$$
où la constante implicite dépend uniquement de $K_1$ et $K_2$.
\end{prop}
\dem
En utilisant la proposition \ref{t131} et la majoration 
\eqref{eq:majH}, 
 on voit sur la définition \eqref{u147} de $\eta(x,v)$ que
\begin{equation}\label{u148}
\eta(x,v) \ll x^{1/2} \int_{0}^{\infty}\frac{|\zeta(\demi+i\tau)|^2}{1+\tau}\min(1,V^{1/2}\big\vert V-\tau\big\vert^{-1})d\tau +x^{1/2} \int_{V}^{\infty}\frac{|\zeta(\demi+i\tau)|^2}{1+\tau^2}d\tau
\end{equation}
avec $V=2\pi x v$. 
Nous utilisons la proposition \ref{t88} pour estimer la première intégrale de \eqref{u148}, et une intégration par parties ainsi que la majoration $I_2(T)\ll T\log(2+T)$ (cf \S\ref{moyquad}) pour estimer la seconde. Nous obtenons la majoration
\[
\eta(x,v) \ll x^{1/2} \frac{\log^2(2+V)}{V^{1/2}}. 
\] 
Si $V < x^2v$, alors 
\[
\eta(x,v) \ll x^{1/2} \frac{\log^2(2+x^2v)}{(xv)^{1/2}}
\ll   \frac{\log^2(2+x^2v)}{\sqrt{x^2v}}. 
\] 
Dans le cas contraire, la décroissance de $u\mapsto \log^2(2+u)/u$ sur $]0,\infty[$ donne 
\[
\eta(x,v) \ll (xV)^{1/2} \frac{\log^2(2+x^2v)}{x^2v}
\ll  \frac{\log^2(2+x^2v)}{\sqrt{x^2v}}.\fine
\]

\medskip

Le dernier résultat de ce paragraphe fournit l'assertion \eqref{t74} de la proposition \ref{prop:eq-approchee-psi1}. 
\begin{prop}\label{prop:pre-final}
 Soit $K_1,K_2>0$. Pour $0<x\ioe K_1$, $v>0$ et $x^2v\soe K_2$, on a
 \[
 \psi(x,v)- x \overline{\psi(1/x,x^2v)}=\demi\log^2\frac 1x+\Big(\gamma-\log 2\pi +i \frac{\pi}{2} \Big)\log \frac 1x+\Fgot(x)+
 O\left(\frac{\log^2(2+x^2v)}{(x^2v)^{1/2}}\right)  
 \] 
 où la constante implicite ne dépend que de $K_1$ et $K_2$. La fonction $\Fgot$ est continue sur $]0,\infty[$ et se prolonge par continuité en $0$.
\end{prop}
 \dem 
 Posons $\kappa = \gamma-\log2\pi+i\pi/2$. D'après la proposition \ref{u140}, on a 
 \begin{align*}
\psi(x,v)- &x \overline{\psi(1/x,x^2v)} \\ 
&=-i\pi\Big(xI(x,v)+\overline{I(1/x,x^2v)}\Big)
+\frac{1}{2} \log^2\frac1x+ \kappa
\log\frac1x 
- \frac{x}{2}\log^2x-
\kappa x \log x
\\ & \quad+ \eps_1(1/x,x^2v)-x\overline{\eps_1(1/x,x^2v)} \pmod \Fcal \, ,
\end{align*}
d'où
\begin{align*}
\psi(x,v)- &x \overline{\psi(1/x,x^2v)} \\ &= -i\pi \eta(x,v)+\frac{1}{2} \log^2\frac1x+ \kappa
\log\frac1x - \frac{x}{2}\log^2x-
\kappa x\log x\\
& \quad+ \eps_1(1/x,x^2v)-x\overline{\eps_1(1/x,x^2v)} +f(x)\expli{où $f\in \Fcal$, d'après la proposition \ref{t100} }\\ 
&= \Fgot(x)+\frac{1}{2} \log^2\frac1x+ \kappa\log\frac1x + O\left(\frac{\log^2(2+x^2v)}{\sqrt{x^2v}}\right),
 \end{align*}
où
$$
\Fgot(x)=f(x)-\frac{x}{2}\log^2x-
\kappa x\log x,
$$ 
et où nous avons utilisé les propositions \ref{u146} et \ref{prop:eta}.\fin

\section{Détermination de la fonction $\Fgot(x)$}\label{t501}

Dans ce paragraphe, la notation $(p.p)$ signifie qu'une identité a lieu pour presque tout nombre réel $x$ de $]0,\infty[$. Nous utiliserons plusieurs fois le fait suivant : si deux fonctions continues sur $]0,\infty[$ coïncident presque partout sur $]0,\infty[$ alors elles sont égales. 
Posons
$$
\psi(x)=\sum_{n\soe 1}\frac{\tau(n)}{n}e^{2\pi inx}.
$$
Comme la série $\sum_{n\soe 1}\tau(n)^2/n^2$ est convergente, la série $\psi(x)$ converge presque partout d'après le théorème de Carleson\footnote{Le résultat découle déjà d'un des premiers théorèmes {\og ancêtres\fg} du théorème de Carleson, dont l'énoncé dû à Jerosch et Weyl (1909, cf. \cite{zbMATH02640541}, p. 78) affirme notamment la convergence presque partout des séries trigonométriques $\sum_n c_ne^{2\pi inx}$ telles que $c_n=O\big((1+\vert n\vert)^{-\kappa}\big)$, avec $\kappa>2/3$.}. L'équation fonctionnelle  de Wilton montre alors que les fonctions $\psi(x)$ et $\widetilde{\psi}(x)=x\overline{\psi(1/x)}$, définies presque partout, vérifient la relation 
$$
\psi(x)-\widetilde{\psi}(x)=\Fgot(x)+\demi\log^2\frac 1x+\Big(\gamma-\log 2\pi +\demi \pi i\Big)\log \frac 1x \quad (p.p.),
$$
relation dont nous pouvons séparer les parties réelle et imaginaire :
\begin{align*}
\psi_1(x)-x\psi_1(1/x)&=\Re\Fgot(x)+\demi\log^2\frac 1x+(\gamma-\log 2\pi )\log \frac 1x \quad (p.p.),\\
\psi_2(x)+x\psi_2(1/x)&=\Im\Fgot(x)+\frac{\pi}{2} \log \frac 1x \quad (p.p.).
\end{align*}

\subsection{Détermination de la fonction $\Im\Fgot(x)$}\label{t110}

Nous proposons deux démonstrations de la relation \eqref{t107} à partir des résultats de \cite{baez-duarte-all}.

\subsubsection{Première démonstration de \eqref{t107}}

Considérons la série
$$
\xi(x)=\sum_{k\soe 1}\frac{B_1(kx)}{k},
$$
où $B_1$ est la première fonction de Bernoulli définie par $B_1(t)=\{t\}-\demi+[t\in\Int]/2$. On sait que cette série converge presque partout et dans $L^2(0,1)$, et a pour somme $-\psi_2/\pi$ (cf. \cite{baez-duarte-all}, p. 222). D'autre part, la proposition 8 de \cite{baez-duarte-all} fournit la relation
$$
\xi(x)+x\xi(1/x)=-A(x)+\frac{1-x}{2}\log x+\frac{x+1}{2}(\log 2\pi-\gamma) \quad (p.p),
$$
d'où \eqref{t107}, en multipliant les deux membres par $-\pi$.

\subsubsection{Deuxième démonstration de \eqref{t107}}\label{t113}

Nous utilisons maintenant la proposition 1 de \cite{baez-duarte-all} : on a l'identité 
\begin{equation}\label{t108}
x\int_x^{\infty}\psi_2(t)\frac{dt}{t^2}=\pi A(x) +\frac{\pi}{2}\log \frac 1x-\frac{\pi}{2}(1-\gamma+\log 2\pi) \quad (x>0).
\end{equation}
En appliquant \eqref{t108} à $1/x$, on obtient
\begin{equation}\label{t109}
\frac 1x\int_0^{x}\psi_2(1/t)dt=\pi A(1/x) +\frac{\pi}{2}\log x-\frac{\pi}{2}(1-\gamma+\log 2\pi) \quad (x>0).
\end{equation}
Comme $A(x)=xA(1/x)$, nous éliminons la fonction $A$ entre \eqref{t108} et \eqref{t109} :
$$
x\int_x^{\infty}\psi_2(t)\frac{dt}{t^2}-\int_0^{x}\psi_2(1/t)dt=\pi\frac{x+1}{2}\log \frac 1x+\pi\frac{x-1}{2}(1-\gamma+\log 2\pi).
$$
En dérivant cette dernière relation, on obtient
$$
\int_x^{\infty}\psi_2(t)\frac{dt}{t^2}-\psi_2(x)/x-\psi_2(1/x)=\frac{\pi}{2}\log \frac 1x-\pi\frac{x+1}{2x}+\frac{\pi}{2}(1-\gamma+\log 2\pi) \quad (p. p.),
$$
d'où
\begin{align*}
\psi_2(x)+x\psi_2(1/x) &=\pi A(x) +\frac{\pi}{2}\log \frac 1x-\frac{\pi}{2}(1-\gamma+\log 2\pi)\\
&\quad -\pi\frac{x}{2}\log \frac 1x + \pi\frac{x+1}{2}-\pi\frac{x}{2}(1-\gamma+\log 2\pi)\\
&=\pi A(x) +\pi\frac{x-1}{2}\log x +\pi\frac{x+1}{2}(\gamma -\log 2\pi)\quad(p.p),
\end{align*}
ce qui achève la démonstration.

\subsection{Détermination de la fonction $\Re\Fgot(x)$}\label{t112}

Notre démarche est d'adapter à $\psi_1$ le raisonnement mené pour $\psi_2$ au \S\ref{t113}. Une observation préalable est que $\psi_1$ est l'opposée de la série trigonométrique conjuguée de $\psi_2$. Cela motive l'emploi de la transformation de Hilbert.

\smallskip

La première étape est d'obtenir une identité analogue à \eqref{t108} pour la fonction $\psi_1$. Commençons avec
$$
\alpha(x)=\int_0^x\psi_1(t)dt=\sum_{n\soe 1}\frac{\tau(n)}{2\pi n^2}\sin 2\pi n x.
$$
Comme $\sin \in \Wcal(-1,0)$ et comme
$$
\sum_{n\soe 1}\frac{\tau(n)}{2\pi n^2}(2\pi n)^{-\sigma}<\infty \quad (\sigma >-1),
$$
la fonction $\alpha$ appartient à $\Wcal(-1,0)$ et
$$
M\alpha(s) =(2\pi)^{-s-1}\zeta^2(s+2)\sin(\pi s/2)\Gamma(s) \quad (-1<\sigma <0).
$$ 
Par conséquent, la fonction $\beta$ définie par
$$
\beta(x)=-\frac{\alpha(x)}{x}+2x\int_x^{\infty}\alpha(t)\frac{dt}{t^3}=x\int_x^{\infty}\psi_1(t)\frac{dt}{t^2}
$$
appartient à $\Wcal(0,1)$ et a pour transformée de Mellin
\begin{align*}
M\beta(s) &=-M\alpha(s-1)+2\frac{M\alpha(s-1)}{s+1}\\
&=\frac{1-s}{1+s}(2\pi)^{-s}\zeta^2(s+1)\sin(\pi (s-1)/2)\Gamma(s-1)\\
&=\frac{\zeta(s+1)}{s(s+1)}\cdot \zeta(s+1) (2\pi)^{-s}\cos (\pi s/2)\Gamma(s+1)\\
&=-\pi\cot(\pi s/2)\cdot\frac{\zeta(s+1)\zeta(-s)}{s(s+1)}\quad (0<\sigma <1),
\end{align*}
où l'on a utilisé l'équation fonctionnelle de $\zeta$ sous la forme \eqref{t78}.
En changeant $s$ en $-s$, on a donc
\begin{equation}\label{t118}
M\beta(-s)=-\pi\cot(\pi s/2)\cdot\frac{\zeta(s)\zeta(1-s)}{s(1-s)}=F(s) \quad (-1<\sigma <0),
\end{equation}
où $F$ est la fonction méromorphe définie par le second membre de \eqref{t125}, et coïncidant avec la transformée de Mellin de la fonction $B$ dans la bande $0<\sigma<1$.

Les deux relations \eqref{t118} et \eqref{t120} permettent d'appliquer la proposition 14 de \cite{baez-duarte-all}. La fonction méromorphe $F$ a un unique pôle sur la droite $\sigma=0$ ; il s'agit d'un pôle triple en $0$ et un calcul standard montre que la partie polaire correspondante est
$$
-\frac{1}{s^3}+\frac{-1+\gamma-\log 2\pi}{s^2}+\frac{c}{s},
$$
où 
$$
c=\frac{\pi^2}{24}-\demi\log^2 2\pi+\demi\gamma^2+\gamma\log 2\pi-\log 2\pi+\gamma +2\gamma_1-1,
$$
$\gamma_1$ désignant la constante de Stieltjes d'indice $1$. En tenant compte de la relation \eqref{t120}, la proposition 14 de \cite{baez-duarte-all} fournit donc
\begin{equation}\label{t121}
B(x)-\beta(1/x)=-\demi\log^2x+(\log 2\pi+1-\gamma)\log x +c \quad (x>0).
\end{equation}

\smallskip

Nous pouvons maintenant reproduire \emph{mutatis mutandis} le raisonnement du \S\ref{t113}. La relation \eqref{t121} appliquée à $x$ et $1/x$ donne
\begin{align}
\int_0^{x}\psi_1(1/t) dt &=xB(x)+\demi x\log^2x+(\log 2\pi+1-\gamma)x\log \frac 1x -cx,\label{t122}\\
x\int_x^{\infty}\psi_1(t)\frac{dt}{t^2}&=B(1/x)+\demi\log^2x+(\log 2\pi+1-\gamma)\log x -c.\label{m123}
\end{align}
Comme 
$$
xB(x)+B(1/x)=-2x\int_0^{\infty}\frac{A(t)}{x+t}\frac{dt}{t}
$$
(cf. \eqref{t133}) est indéfiniment dérivable sur $]0,\infty[$, on peut dériver l'identité obtenue en ajoutant \eqref{t122} et \eqref{m123} et obtenir ainsi pour presque tout $x>0$,
\begin{align*}
\psi_1(1/x)-\psi_1(x)/x+\int_x^{\infty}\psi_1(t)\frac{dt}{t^2}&=-2\int_0^{\infty}\frac{A(t)}{(x+t)^2}dt+
\frac{\log x +\log 2\pi+1-\gamma}{x} \\ &\quad +\demi\log^2 x+(\log 2\pi -\gamma)\log \frac 1x -\log 2\pi -1 +\gamma-c 
\end{align*}
puis, en réutilisant \eqref{m123} :
\begin{align*}
x\psi_1(1/x)-\psi_1(x)&=-2x\int_0^{\infty}\frac{A(t)}{(x+t)^2}dt-B(1/x)-\demi\log^2x-(\log 2\pi-\gamma)\log x \\
&\quad+x\Big(\demi \log^2 x+(\log 2\pi -\gamma)\log \frac 1x -\log 2\pi -1 +\gamma-c\Big)\\&\quad+c+\log 2\pi+1-\gamma\quad(p.p).
\end{align*}
On obtient bien \eqref{t111} en utilisant l'identité \eqref{t134} :
$$
2x\int_0^{\infty}\frac{A(t)}{(x+t)^2}dt+B(1/x)=2x\cdot\vp\int_0^{\infty}A(t)\frac{x^2+t^2}{(t-x)(t+x)^2}\, \frac{dt}{t}
$$
et la relation
$$
c+\log 2\pi+1-\gamma=\frac{\pi^2}{24}-\demi\ln^2 2\pi+\demi\gamma^2+\gamma\ln 2\pi+2\gamma_1.
$$

\begin{center}
  {\sc Remerciements}
\end{center}
\begin{quote}
{\footnotesize Outre leurs laboratoires respectifs, les auteurs remercient les institutions ayant favorisé leur travail sur cet article : le laboratoire franco-russe Poncelet (CNRS, Université Indépendante de Moscou), et l'Institut Mittag-Leffler (Djursholm).

Les auteurs ont également le plaisir de remercier l'arbitre anonyme pour sa lecture attentive de cet article, ses utiles suggestions, et pour le contenu du \S\ref{t137}.
}
\end{quote}

\medskip


\begin{thebibliography}{10}

\bibitem{baez-duarte-all}
{\scshape L.~B{\'a}ez-Duarte, M.~Balazard, B.~Landreau {\normalfont
  \smfandname} E.~Saias} -- {\og \'{E}tude de l'autocorr\'elation
  multiplicative de la fonction {\og partie fractionnaire\fg}\fg},
  \emph{Ramanujan J.} \textbf{9} (2005), p.~215--240.

\bibitem{bourgain-2014}
{\scshape J.~Bourgain} -- {\og Decoupling, exponential sums and the {R}iemann
  {Z}eta function\fg}, \emph{arXiv:1408.5794v1} (2014).

\bibitem{chandrasekharan-70}
{\scshape K.~Chandrasekharan} -- \emph{Arithmetical functions}, Die Grundlehren
  der mathematischen Wissenschaften, vol. 167, Springer, Berlin, Heidelberg,
  New York, 1970.

\bibitem{HL}
{\scshape G.~H. Hardy {\normalfont \smfandname} J.~E. Littlewood} -- {\og Some
  problems of diophantine approximation\fg}, \emph{Acta Math.} \textbf{37}
  (1914), p.~193--239, II The trigonometrical series associated with the
  elliptic $\vartheta$-function.

\bibitem{huxley-2005}
{\scshape M.~N. Huxley} -- {\og Exponential sums and the {R}iemann {Z}eta
  function. {V}\fg}, \emph{Proc. London Math. Soc. (3)} \textbf{90} (2005),
  p.~1--41.

\bibitem{53.0313.01}
{\scshape A.~E. Ingham} -- {\og {Mean-value theorems in the theory of the
  Riemann Zeta-function}\fg}, \emph{Proc. London Math. Soc. (2)} \textbf{27}
  (1927), p.~273--300.

\bibitem{ivic}
{\scshape A.~Ivi{\'c}} -- \emph{The {R}iemann zeta-function}, Dover
  Publications Inc., Mineola, NY, 2003.

\bibitem{zbMATH02640541}
{\scshape F.~{Jerosch} {\normalfont \smfandname} H.~{Weyl}} -- {\og {\"Uber die
  Konvergenz von Reihen, die nach periodischen Funktionen Fortschreiten.}\fg},
  \emph{{Math. Ann.}} \textbf{66} (1909), p.~67--80.

\bibitem{king-09}
{\scshape F.~W. King} -- \emph{Hilbert transforms, vol. 1}, Encyclopedia of
  mathematics and its applications, vol. 124, Cambridge University Press,
  Cambridge, 2009.

\bibitem{37.0450.01}
{\scshape N.~Nielsen} -- \emph{{Handbuch der Theorie der Gammafunktionen.}},
  {B. G. Teubner, Leipzig}, 1906.

\bibitem{rivoal-2012}
{\scshape T.~Rivoal} -- {\og On the convergence of diophantine {D}irichlet
  series\fg}, \emph{Proc. Edimb. Math. Soc.} \textbf{55} (2012), p.~513--541.


\bibitem{rivoal-roques-2013}
{\scshape T.~Rivoal {\normalfont \smfandname} J.~Roques} -- {\og Convergence
  and modular type properties of a twisted {R}iemann series\fg}, \emph{Unif.
  Distrib. Theory} \textbf{8} (2013), p.~97--119.

\bibitem{rivoal-seuret-2014}
{\scshape T.~Rivoal {\normalfont \smfandname} S.~Seuret} -- {\og
  {H}ardy-{L}ittlewood series and even continued fractions\fg}, \emph{J. Anal.
  Math.} \textbf{125} (2015), p.~175--225.



\bibitem{MR1366197}
{\scshape G.~Tenenbaum} -- \emph{Introduction \`a la th\'eorie analytique et
  probabiliste des nombres}, 3\up{e} \smfedname, Belin, Paris, 2008.

\bibitem{titchmarsh-fourier}
{\scshape E.~C. Titchmarsh} -- \emph{Introduction to the theory of {F}ourier
  integrals}, second \smfedname, Clarendon Press, Oxford, 1948.

\bibitem{titchmarsh-zeta}
\bysame , \emph{The theory of the {R}iemann zeta function}, second \smfedname,
  Clarendon Press, Oxford, 1986, revised by {D. R. H}eath-{B}rown.

\bibitem{34.0231.03}
{\scshape G.~Vorono\"i} -- {\og {Sur un probl\`eme du calcul des fonctions
  asymptotiques}\fg}, \emph{J. reine angew. Math} \textbf{126} (1903),
  p.~241--282.

\bibitem{voronoi1}
{\scshape G.~Vorono{\"{\i}}} -- {\og Sur une fonction transcendante et ses
  applications \`a la sommation de quelques s\'eries\fg}, \emph{Ann. Sci.
  \'Ecole Norm. Sup. (3)} \textbf{21} (1904), p.~207--267, 459--533.

\bibitem{Wilton}
{\scshape J.~Wilton} -- {\og An approximate functional equation with
  applications to a problem of diophantine approximation\fg}, \emph{J. reine
  angew. Math} \textbf{169} (1933), p.~219--237.

\bibitem{MR1963498}
{\scshape A.~Zygmund} -- \emph{Trigonometric series. {V}ol. {I}, {II}}, third
  \smfedname, Cambridge Mathematical Library, Cambridge University Press,
  Cambridge, 2002.

\end{thebibliography}
\providecommand{\bysame}{\leavevmode ---\ }
\providecommand{\og}{``}
\providecommand{\fg}{''}
\providecommand{\smfandname}{et}
\providecommand{\smfedsname}{\'eds.}
\providecommand{\smfedname}{\'ed.}
\providecommand{\smfmastersthesisname}{M\'emoire}
\providecommand{\smfphdthesisname}{Th\`ese}

\begin{multicols}{2}
\footnotesize

\noindent BALAZARD, Michel\\
Institut de Math\'ematiques de Marseille\\
CNRS, Universit\'e d'Aix-Marseille\\
Campus de Luminy, Case 907\\
13288 Marseille Cedex 9\\
FRANCE\\
Adresse \'electronique : \texttt{balazard@math.cnrs.fr}

\smallskip

\noindent MARTIN, Bruno\\
Laboratoire de Math\'ematiques Pures et Appliqu\'ees\\
CNRS, Universit\'e du Littoral C\^ote d'Opale\\
50 rue F. Buisson, BP 599\\
62228 Calais Cedex\\
FRANCE\\
Adresse \'electronique : \texttt{martin@lmpa.univ-littoral.fr}
\end{multicols}

\end{document}